\documentclass[12pt]{article}
\usepackage{oldgerm}
\usepackage{euler}
\usepackage{graphics}
\usepackage{bm}
\usepackage{graphicx}
\usepackage{amsmath}
\usepackage{amssymb}
\usepackage{amstext}
\usepackage{amscd}
\usepackage{amsfonts}
\usepackage{appendix}
\usepackage{color}
\DeclareMathAlphabet{\EuFrak}{U}{euf}{m}{n}
\DeclareMathAlphabet{\EuScript}{U}{eus}{m}{n}

\newcommand{\nd}{\noindent}

\newcommand{\be}{\begin{equation}}
\newcommand{\ee}{\end{equation}}
\newcommand{\ben}{\begin{eqnarray}}
\newcommand{\een}{\end{eqnarray}}

\title{{\bf q-Gamow States as continuous linear functionals
on analytical test functions}}

\author{{A. Plastino$^1$, M.C.Rocca$^1$}\\
\small{$^1$ Department of Physics,
La Plata National University}\\ \small{ and
Argentina's National Research Council}\\
\small{(IFLP-CCT-CONICET)-C. C. 727, 1900 La Plata - Argentina}}

\date{\today}

\begin{document}

\maketitle

\begin{abstract}

We define here q-Gamow states corresponding to Tsallis'
q-statistics. We compute for them their norm, mean energy value an
the q-analogue of the  Breit-Wigner distribution (a
q-Breit-Wigner).

\nd {\bf Keywords:} Gamow States, q-Gamow States,
Ultradistributions.

\end{abstract}

\newpage

\renewcommand{\theequation}{\arabic{section}.\arabic{equation}}

\section{Introduction}

In three previous papers \cite{tp1,tp2,tp3} we have shown that
Gamow-states  \cite{tp4} can be interpreted as Sebastiao e Silva's
Ultradistributions  \cite{tp5,tp6,tp7}, whose proper treatment
appeals to  Rigged Hilbert Space \cite{tp8,tp9,tp10}.

Lately, one finds many high energy experiments that can be
interpreted via Tsallis' q-statistics  \cite{tp11o}. Indeed, there
has been  increased use by  LHC experiments of such q-statistics
and, specially, of the distribution associated with a stationary
state within q-statistics, that seems to describe very well the
transverse momentum distributions of all different types of
hadrons. All four LHC experiments have published results for these
distributions that are well fitted by the q-exponential function.
The resulting value of q is around 1.15, quite different from
Shannon's-Boltzmann's  $q=1$. This means that the stationary
states of the particles before the hadronization are not in
thermal equilibrium. Moreover, the distribution is very robust and
practically the same for different hadrons,  spanning a range of
different energies. Maybe, one of the most impressive results,
recently published, is the measurement of the $p_T$ distribution
over a logarithmic range of 14 decades. It was found that the same
expression of a q exponential ($q=1.15$) fits the data over the
full range of these fourteen decades. A theory that fits a range
of couple of decades is already very interesting but doing so with
such a large range of decades, with the same distribution, is rare
indeed  (see, for instance, \cite{tp11,phenix}).

 These circumstances strongly motivate us to investigate complex energy states
  related to
  the  q-exponential distributions, that is,  q-Gamow
states, and establish their relation with Gamow-states.
 Er focus attention then on decay states at a great distance from the dispersion center and ascertain
 that a q-Gamow representation is adequate.

\setcounter{equation}{0}

\section{Gamow States}

Following  \cite{tp1,tp2,tp3} we define a Gamow-state in a
dispersion-less space  as
\begin{equation}
\label{eq2.1} |\psi_G>=\int\limits_{-\infty}^{\infty} \left\{{\cal
H}[\Im(p)]{\cal H}(x)-{\cal H}[-\Im(p)]{\cal H}(-x)\right\}
e^{\frac {ipx} {\hbar}}|x>\;dx,
\end{equation}
or
\begin{equation}
\label{eq2.2} \psi_G(x)= \left\{{\cal H}[\Im(p)]{\cal H}(x)-{\cal
H}[-\Im(p)]{\cal H}(-x)\right\} e^{\frac {ipx} {\hbar}}.
\end{equation}
The norm-squared for such a state reads
\begin{equation}
\label{eq2.3} <\psi_G|\psi_G>=\int\limits_0^{\infty} {\cal
H}[\Im(p)]e^{\frac {i(p-p^{\ast})x} {\hbar}}\;dx -
\int\limits_{-\infty}^0 {\cal H}[-\Im(p)]e^{\frac {i(p-p^{\ast})x}
{\hbar}}\;dx.
\end{equation}
These  integrals can be easily evaluated. One finds
\begin{equation}
\label{eq2.4} <\psi_G|\psi_G>=\left\{{\cal H}[\Im(p)]- {\cal
H}[-\Im(p)]\right\}\frac {\hbar} {i(p^{\ast}-p)}= \frac {\hbar}
{2|\Im(p)|}.
\end{equation}
Accordingly, the normalized Gamow-state  $\phi_G$ becomes
\begin{equation}
\label{eq2.5} |\phi_G>=\sqrt{\frac {2|\Im(p)|} {\hbar}}|\psi_G>.
\end{equation}
Since
\begin{equation}
\label{eq2.6} <\phi_G|(H|\phi_G>)=\frac {p^2} {2m},
\end{equation}
\begin{equation}
\label{eq2.7} (<\phi_G|H)|\phi_G>=\frac {p^{\ast 2}} {2m},
\end{equation}
one encounters, for the energy mean value \cite{tp1,tp2,tp3}
\begin{equation}
\label{eq2.8} <H>=\frac {1}
{2}\left[<\phi_G|(H|\phi_G>)+(<\phi_G|H)|\phi_G>\right]= \frac
{p^2+p^{\ast 2}} {4m}= \frac {\Re(p^2)} {2m}.
\end{equation}

%%%%%%%para

in order to obtain de probability distribution associated to a
q-Gamow state we start by the looking at scalar product between
this state and a free one:
\begin{equation}
\label{eq2.9} <\phi|\phi_G>=\frac {1} {\hbar}\sqrt{\frac
{|\Im(p)|} {\pi}} \left\{\int\limits_0^{\infty} {\cal
H}[\Im(p)]e^{\frac {i(p-k)x} {\hbar}}\;dx- \int\limits_{-\infty}^0
{\cal H}[-\Im(p)]e^{\frac {i(p-k)x} {\hbar}}\;dx\right\}.
\end{equation}
Thus,
\begin{equation}
\label{eq2.10}
<\phi|\phi_G>=\frac {i\sqrt{\frac {|\Im(p)|} {\pi}}}
{p-k}
\end{equation}
The ensuing probability distribution is the  Breit-Wigner one
\begin{equation}
\label{eq2.11} |<\phi|\phi_G>|^2=\frac {|\Im(p)|}
{\pi\left\{[\Re(p)-k]^2+\Im(p)^2\right\}}.
\end{equation}

\setcounter{equation}{0}

\section{q-Gamow States}

According to the q-statistics strictures (this word exits!)
 we must replace everywhere ordinary exponentials by so.called q-exponentials $e_q(x)$ \cite{tp11o}
 \be e_q(x) = [1+ (1-q)x]^{1/1-q}; \,\,q \in \mathcal{R}, \ee
that becomes the ordinary exponential at $q=1$. Accordingly,

\[|\psi_{qG}>=\int\limits_{-\infty}^{\infty}
\left\{{\cal H}[\Im(p)]{\cal H}(x)-{\cal H}[-\Im(p)]{\cal H}(-x)\right\}
\otimes\]
\begin{equation}
\label{eq3.1} \left[1-\frac {i(q-1)px}
{\hbar\sqrt{2(q+1)}}\right]^{\frac {2} {1-q}} |x>\;dx,
\end{equation}
or

\begin{equation}
\label{eq3.2} \psi_{qG}(x)= \left\{{\cal H}[\Im(p)]{\cal
H}(x)-{\cal H}[-\Im(p)]{\cal H}(-x)\right\} \left[1-\frac
{i(q-1)px} {\hbar\sqrt{2(q+1)}}\right]^{\frac {2} {1-q}}.
\end{equation}
The  norm of a  q-Gamow state is

\[<\psi_{qG}|\psi_{qG}>=\int\limits_0^{\infty}
{\cal H}[\Im(p)]
\left[1-\frac {i(q-1)px} {\hbar\sqrt{2(q+1)}}\right]^{\frac {2} {1-q}}
\left[1+\frac {i(q-1)p^{\ast}x} {\hbar\sqrt{2(q+1)}}\right]^{\frac {2} {1-q}}
\;dx \]
\begin{equation}
\label{eq3.3} +\int\limits_{-\infty}^0 {\cal H}[-\Im(p)]
\left[1-\frac {i(q-1)px} {\hbar\sqrt{2(q+1)}}\right]^{\frac {2}
{1-q}} \left[1+\frac {i(q-1)p^{\ast}x}
{\hbar\sqrt{2(q+1)}}\right]^{\frac {2} {1-q}} \;dx,
\end{equation}
or equivalently,
\[<\psi_{qG}|\psi_{qG}>=\int\limits_0^{\infty}
{\cal H}[\Im(p)]
\left[1+\frac {2(q-1)\Im(p)x} {\hbar\sqrt{2(q+1)}}+
\frac {(q-1)^2|p|^2x^2} {\hbar^22(q+1)}
\right]^{\frac {2} {1-q}}\;dx\]
\begin{equation}
\label{eq3.4} +\int\limits_0^{\infty} {\cal H}[-\Im(p)]
\left[1-\frac {2(q-1)\Im(p)x} {\hbar\sqrt{2(q+1)}}+ \frac
{(q-1)^2|p|^2x^2} {\hbar^22(q+1)} \right]^{\frac {2} {1-q}}\;dx,
\end{equation}
that can be recast as
\begin{equation}
\label{eq3.5} <\psi_{qG}|\psi_{qG}>= \int\limits_0^{\infty}
\left[1+\frac {2(q-1)|\Im(p)|x} {\hbar\sqrt{2(q+1)}}+ \frac
{(q-1)^2|p|^2x^2} {\hbar^22(q+1)} \right]^{\frac {2} {1-q}}\;dx.
\end{equation}
We effect now the change of variables $y=\frac {(q-1)|p|x} {\hbar
\sqrt{2(q+1)}}$ and obtain
\begin{equation}
\label{eq3.6} <\psi_{qG}|\psi_{qG}>= \frac {\hbar\sqrt{2(q+1)}}
{(q-1)|p|} \int\limits_0^{\infty} \left[1+\frac {2|\Im(p)|y}
{|p|^2}+y^2 \right]^{\frac {2} {1-q}}\;dy.
\end{equation}
After a new change of variables  $z=y+\frac {|\Im(p)|} {|p|}$
 we find
\begin{equation}
\label{eq3.7} <\psi_{qG}|\psi_{qG}>= \frac {\hbar\sqrt{2(q+1)}}
{(q-1)|p|} \int\limits_{\frac {|\Im(p)|} {|p|}}^{\infty}
\left\{z^2+\frac {[\Re(p)]^2} {|p|^2} \right\}^{\frac {2}
{1-q}}\;dz.
\end{equation}
Finally, after a third change of variables  $s=z^2$ we get for our
norm
\begin{equation}
\label{eq3.8} <\psi_{qG}|\psi_{qG}>= \frac {\hbar\sqrt{2(q+1)}}
{2(q-1)|p|} \int\limits_{\frac {|\Im(p)|^2} {|p|^2}}^{\infty}
s^{-\frac {1} {2}} \left\{s+\frac {[\Re(p)]^2} {|p|^2}
\right\}^{\frac {2} {1-q}}\;ds.
\end{equation}
Using the result given in \cite{tp12} we arrive to:
\[<\psi_{qG}|\psi_{qG}>=
\frac {\hbar} {5-q}
\frac {\sqrt{2(q+1)}} {|p|}
\left\{\frac {[\Im(p)]^2} {|p|^2}\right\}^{\frac {q-5} {2(q-1)}}\otimes\]
\begin{equation}
\label{eq3.9} F\left(\frac {2} {q-1},\frac {5-q} {2(q-1)}; \frac
{3+q} {2(q-1)};-\frac {[\Re(p)]^2} {[\Im(p)]^2}\right).
\end{equation}
It is shown in  \cite{tp13} that
\[F\left(\frac {2} {q-1},\frac {5-q} {2(q-1)};
\frac {3+q} {2(q-1)};-\frac {[\Re(p)]^2} {[\Im(p)]^2}\right)=\]
\begin{equation}
\label{eq3.10} \left\{\frac {|p|^2} {[\Im(p)]^2}\right\}^{\frac
{q-5} {2(q-1)}} F\left(\frac {1} {2},\frac {5-q} {2(q-1)}; \frac
{3+q} {2(q-1)};\frac {[\Re(p)]^2} {|p|^2}\right),
\end{equation}
which yields for the norm the expression
\[<\psi_{qG}|\psi_{qG}>=
\frac {\hbar} {5-q}
\frac {\sqrt{2(q+1)}} {|p|}\otimes\]
\begin{equation}
\label{eq3.11} F\left(\frac {1} {2},\frac {5-q} {2(q-1)}; \frac
{3+q} {2(q-1)};\frac {[\Re(p)]^2} {|p|^2}\right)= [A(q,p)]^2,
\end{equation}
so that the normalized  q-Gamow state becomes
\begin{equation}
\label{eq3.12} |\phi_{qG}>=[A(q,p)]^{-1}|\psi_{qG}>.
\end{equation}
Noticing that
\begin{equation}
\label{eq3.13}
\lim_{q\rightarrow 1}
F\left(\frac {1} {2},\frac {5-q} {2(q-1)};
\frac {3+q} {2(q-1)};\frac {[\Re(p)]^2} {|p|^2}\right)=
F\left(\frac {1} {2},4;4;\frac {[\Re(p)]^2} {|p|^2}\right)
\end{equation}
and using a result of \cite{tp14} one has
\begin{equation}
\label{eq3.14} F\left(\frac {1} {2},4;4;\frac {[\Re(p)]^2}
{|p|^2}\right)= \left[\frac {[\Im(p)]^2} {|p|^2}\right]^{-\frac
{1} {2}},
\end{equation}
and
\begin{equation}
\label{eq3.15} \lim_{q\rightarrow 1}[A(q,p)]^2=\frac {\hbar}
{2|\Im(p)|}.
\end{equation}
Using now, from  \cite{tp15},
\begin{equation}
\label{eq3.16} H\phi_q(x)=\frac {p^2} {2m}[\phi_q(x)]^q,
\end{equation}
we encounter
\[<\phi_{qG}|(H|\phi_{qG}>)=[A(p,q)]^{-2}\frac {p^2} {2m}\otimes\]
\[\left\{
\int\limits_0^{\infty} {\cal H}[\Im(p)]
\left[1-\frac {i(q-1)px} {\hbar\sqrt{2(q+1)}}\right]^{\frac {2q} {1-q}}
\left[1+\frac {i(q-1)p^{\ast}x} {\hbar\sqrt{2(q+1)}}\right]^{\frac {2} {1-q}}
\;dx \right.\]
\begin{equation}
\label{eq3.17} \left.+\int\limits_{-\infty}^0 {\cal H}[-\Im(p)]
\left[1-\frac {i(q-1)px} {\hbar\sqrt{2(q+1)}}\right]^{\frac {2q}
{1-q}} \left[1+\frac {i(q-1)p^{\ast}x}
{\hbar\sqrt{2(q+1)}}\right]^{\frac {2} {1-q}} \;dx\right\},
\end{equation}
or equivalently,
\[<\phi_{qG}|(H|\phi_{qG}>)=[A(p,q)]^{-2}\frac {p^2} {2m}\otimes\]
\[\left\{
\int\limits_0^{\infty} {\cal H}[\Im(p)]
\left[1-\frac {i(q-1)px} {\hbar\sqrt{2(q+1)}}\right]^{\frac {2q} {1-q}}
\left[1+\frac {i(q-1)p^{\ast}x} {\hbar\sqrt{2(q+1)}}\right]^{\frac {2} {1-q}}
\;dx \right.\]
\begin{equation}
\label{eq3.18} \left.+\int\limits_0^{\infty} {\cal H}[-\Im(p)]
\left[1+\frac {i(q-1)px} {\hbar\sqrt{2(q+1)}}\right]^{\frac {2q}
{1-q}} \left[1-\frac {i(q-1)p^{\ast}x}
{\hbar\sqrt{2(q+1)}}\right]^{\frac {2} {1-q}} \;dx\right\}.
\end{equation}
We can recast  (\ref{eq3.18}) as
\[<\phi_{qG}|(H|\phi_{qG}>)=[A(p,q)]^{-2}\frac {p^2} {2m}\otimes\]
\[\left\{{\cal H}[\Im(p)]
\left[\frac {i(q-1)p^{\ast}} {\hbar\sqrt{2(q+1)}}\right]^{
\frac {2} {1-q}}
\left[\frac {-i(q-1)p} {\hbar\sqrt{2(q+1)}}\right]^{
\frac {2q} {1-q}}\right.\otimes\]
\[\int\limits_0^{\infty}
\left[1-\frac {i(q-1)px} {\hbar\sqrt{2(q+1)}}\right]^{\frac {2q} {1-q}}
\left[1+\frac {i(q-1)p^{\ast}x} {\hbar\sqrt{2(q+1)}}\right]^{\frac {2} {1-q}}
\;dx \]
\[+{\cal H}[-\Im(p)]
\left[\frac {-i(q-1)p^{\ast}} {\hbar\sqrt{2(q+1)}}\right]^{
\frac {2} {1-q}}
\left[\frac {i(q-1)p} {\hbar\sqrt{2(q+1)}}\right]^{
\frac {2q} {1-q}}\otimes\]
\begin{equation}
\label{eq3.19} \left.\int\limits_0^{\infty} \left[1+\frac
{i(q-1)px} {\hbar\sqrt{2(q+1)}}\right]^{\frac {2q} {1-q}}
\left[1-\frac {i(q-1)p^{\ast}x} {\hbar\sqrt{2(q+1)}}\right]^{\frac
{2} {1-q}} \;dx\right\}.
\end{equation}
We use now a result from  \cite{tp16} and obtain
\[<\phi_{qG}|(H|\phi_{qG}>)=-\frac {\hbar} {i[A(p,q)]^2}
\frac {p} {2m}{\cal B}\left(1,\frac {3+q} {q-1}\right)
\frac {\sqrt{2(q+1)}} {(q-1)}
\otimes\]
\begin{equation}
\label{3.20} \left\{{\cal H}[\Im(p)]-{\cal H}[-\Im(p)]\right\}
F\left(1,\frac {2} {q-1};\frac {2(q+1)} {q-1}; 1+\frac {p^{\ast}}
{p}\right),
\end{equation}
or equivalently,
\[<\phi_{qG}|(H|\phi_{qG}>)=-\frac {\hbar}
{i[A(p,q)]^2} \frac {p} {2m} \frac {\sqrt{2(q+1)}} {(3+q)}
\otimes\]
\begin{equation}
\label{3.21} \left\{{\cal H}[\Im(p)]-{\cal H}[-\Im(p)]\right\}
F\left(1,\frac {2} {q-1};\frac {2(q+1)} {q-1}; 1+\frac {p^{\ast}}
{p}\right).
\end{equation}
In analogous fashion we find
\[(<\phi_{qG}|H)|\phi_{qG}>=\frac {\hbar} {i[A(p,q)]^2}
\frac {p^{\ast}} {2m}
\frac {\sqrt{2(q+1)}} {(3+q)}
\otimes\]
\begin{equation}
\label{3.22} \left\{{\cal H}[\Im(p)]-{\cal H}[-\Im(p)]\right\}
F\left(1,\frac {2} {q-1};\frac {2(q+1)} {q-1}; 1+\frac {p}
{p^{\ast}}\right).
\end{equation}
Thus, according to  \cite{tp1,tp2,tp3} we obtain for the mean
energy value
\begin{equation}
\label{3.23} <H>_q=\frac {1} {2} \left[<\phi_{qG}|(H|\phi_{qG}>)+
(<\phi_{qG}|H)|\phi_{qG}>\right].
\end{equation}
Additionally, since
\begin{equation}
\label{3.24} \lim_{q\rightarrow 1} F\left(1,\frac {2} {q-1};\frac
{2(q+1)} {q-1}; 1+\frac {p^{\ast}} {p}\right)= \frac {2p}
{p-p^{\ast}},
\end{equation}
we have
\begin{equation}
\label{3.25} \lim_{q\rightarrow 1}<H>_q=\frac {\Re(p^2)} {2m}=<H>.
\end{equation}
\vskip 3mm

We investigate now the q-analogue of the Breit-Wigner distribution
tackling
\[<\phi|\phi_{Gq}>=\frac {1} {\sqrt{2\pi\hbar A(q,p)}}\left\{
{\cal H}[\Im(p)]\int\limits_0^{\infty}e^{-ikx}
\left[1+\frac {i(1-q)px} {\hbar\sqrt{2(q+1)}}\right]^{
\frac {2} {1-q}}\;dx\right.\]
\begin{equation}
\label{eq3.26} -\left.{\cal
H}[-\Im(p)]\int\limits_0^{\infty}e^{ikx} \left[1-\frac {i(1-q)px}
{\hbar\sqrt{2(q+1)}}\right]^{ \frac {2} {1-q}}\;dx\right\},
\end{equation}
and rewrite it as
\[<\phi|\phi_{Gq}>=\frac {1} {\sqrt{2\pi\hbar A(q,p)}}\left\{
H[\Im(p)]\left[\frac {i(1-q)p} {\hbar\sqrt{2(q+1)}}\right]^{
\frac {2} {1-q}}\right.\otimes\]
\[\int\limits_0^{\infty}e^{-ikx}
\left[x+\frac {\hbar\sqrt{2(q+1)}} {i(1-q)p}\right]^{
\frac {2} {1-q}}\;dx-
H[-\Im(p)]\left[-\frac {i(1-q)p} {\hbar\sqrt{2(q+1)}}\right]^{
\frac {2} {1-q}}\otimes\]
\begin{equation}
\label{eq3.27} \left.\int\limits_0^{\infty}e^{ikx} \left[x-\frac
{\hbar\sqrt{2(q+1)}} {i(1-q)p}\right]^{ \frac {2}
{1-q}}\;dx\right\}.
\end{equation}
We appeal now to a result of  \cite{tp17} and obtain
\[<\phi|\phi_{Gq}>=-i\sqrt{\frac {\hbar} {2\pi A(q,p)}}
\left[\frac {\sqrt{2(q+1)}} {(1-q)p}\right]^{
\frac {2} {q-1}}k^{\frac {3-q} {q-1}}
e^{\frac {\sqrt{2(q+1)}k} {(1-q)p}}\otimes\]
\begin{equation}
\label{eq3.28} \Gamma\left[\frac {3-q} {1-q}, \frac
{\sqrt{2(q+1)}\;k} {(1-q)p}\right],
\end{equation}
which leads to the q-Breit-Wigner result
\[|<\phi|\phi_{Gq}>|^2=\frac {\hbar} {2\pi A(q,p)}
\left[\frac {2(q+1)} {(1-q)^2|p|^2}\right]^{
\frac {2} {q-1}}k^{\frac {2(3-q)} {q-1}}
e^{\frac {\sqrt{2(q+1)}k(p+p^{\ast})} {(1-q)|p|^2}}\otimes\]
\begin{equation}
\label{eq3.29} \Gamma\left[\frac {3-q} {1-q}, \frac
{\sqrt{2(q+1}\;k} {(1-q)p}\right] \left\{\Gamma\left[\frac {3-q}
{1-q}, \frac {\sqrt{2(q+1}\;k} {(1-q)p}\right]\right\}^{\ast}.
\end{equation}
Note that (\ref{eq3.26}) converges uniformly for $q\rightarrow 1$
[since the q-exponential converges in that way to the ordinary
one], i.e.,
\begin{equation}
\label{eq3.30} \lim_{q\rightarrow 1}|<\phi|\phi_{Gq}>|^2=\frac
{|\Im(p)|} {\pi\left\{[\Re(p)-k]^2+\Im(p)^2\right\}}.
\end{equation}

\setcounter{equation}{0}

\section{Conclusions}

In this work we gave introduced q-Gamow states. For that purpose
we have computed their norm, the mean energy value, and the
concomitant   q-Breit-Wigner distributions.
 In all instance, results tend to the customary ones for  $q\rightarrow 1$.

\end{document}